\documentclass{article}
\author{Claire Levaillant\\clairelevaillant@yahoo.fr}
\title{Classification of the invariant subspaces of the Cohen--Wales
representation of the Artin group of type $D_n$}

\usepackage{amsmath, amssymb}
\usepackage{amsthm}
\usepackage{palatino}
\usepackage{amscd}
\usepackage{graphicx}
\usepackage{epsfig}
\usepackage{verbatim}

\newcommand{\n}{\nu}
\newcommand{\Q}{\mathbb{Q}}

\newcommand{\la}{\lambda}

\newcommand{\ovl}{\overline}
\newcommand{\unsurr}{\frac{1}{r}}
\newcommand{\unsur}{\frac{1}}

\newcommand{\T}{\mathcal{T}}
\newcommand{\U}{\mathcal{U}}
\newcommand{\V}{\mathcal{V}}
\newcommand{\IH}{\mathcal{H}}

\newcommand{\lb}{\lbrace}
\newcommand{\rb}{\rbrace}
\newcommand{\noin}{\noindent}

\newcommand{\ih}{\IH_{F,r^2}}

\newcommand{\W}{\mathcal{W}}

\newcommand{\da}{\downarrow}

\newcommand{\B}{\mathcal{B}}

\newcommand{\g}{\gamma}

\newcommand{\cil}{\frac{n(n-3)}{2}}
\newcommand{\dbw}{\frac{(n-1)(n-2)}{2}}
\newcommand{\chl}{\frac{n(n-1)}{2}}

\newcommand{\cgwn}{\mathcal{C}\mathcal{G}\mathcal{W}(D_n)}
\newcommand{\wh}{\widehat}
\newcommand{\mcalh}{\mathcal{H}}

\newcommand{\h}{(\mcalh)_}

\begin{document}
\maketitle \noin\textbf{Abstract.} Recently, Cohen and Wales built a
faithful linear representation of the Artin group of type $D_n$,
hence showing the linearity of this group. It was later discovered
that this representation is reducible for some complex values of its
two parameters. It was also shown that when the representation is
reducible, the action on a proper invariant subspace is a Hecke
algebra action of type $D_n$. The goal of this paper is to classify
these proper invariant subspaces in terms of Specht modules indexed
by double partitions of the integer $n$. This work is the
continuation of \cite{CL}.

\section{Introduction}
In \cite{CW}, Cohen and Wales built representations of the Artin
groups of types $A$, $D$ and $E$. In type $A$, their representation
is equivalent to the Lawrence--Krammer representation of the braid
group, a famous representation that was used in \cite{BIG} and
independently in \cite{KR} to show the linearity of the braid group.
In \cite{CL}, we use knot theory to construct a representation
$\n^{(n)}$ of an algebra that contains the Artin group of type $D_n$
and that depends on two parameters $l$ and $m$. This algebra is
defined by the authors in \cite{CGW} and is a generalization of the
Birman-Murakami-Wenzl algebra to type $D_n$. We thus call it the CGW
algebra of type $D_n$. We show that as a representation of the Artin
group of type $D_n$ and up to some change of parameters, this
representation is equivalent to the Cohen-Wales representation of
the Artin group type $D_n$. We prove that the representation is
generically irreducible, but that when the parameters are specified
to some nonzero complex numbers, it becomes reducible. We give a
reducibility criterion for this representation, thus obtaining a
reducibility criterion for the original representation of Cohen and
Wales. The parameters $t$ and $r$ of the Cohen-Wales representation
are related to the parameters $l$ and $m$ of the algebra by
$m=r-\unsurr$ and $l=\unsur{t\,r^3}$. The complex parameters $t$ and
$r$ for which the Cohen-Wales representation is reducible are given
in the Main Theorem of \cite{CL}. It is also shown in \cite{CL} that
when the representation is reducible the action on a proper
invariant subspace is a Hecke algebra action. The $r$ of \cite{CL}
and of the current paper is the $\unsurr$ of the Cohen-Wales
representation. As in \cite{CL}, we denote by $\mcalh(D_n)$ the
Hecke algebra of type $D_n$ with parameter $r^2$ over the field
$\Q(l,r)$. The classes of irreducible $\mcalh(D_n)$-modules are
called Specht modules and are indexed by double partitions of the
integer $n$, see $\S\,3.2$ of \cite{CL}. We denote by $\ih(n)$ the
Iwahori-Hecke agebra of the symmetric group $Sym(n)$ with parameter
$r^2$ over the field $\Q(l,r)$. Some important aspects of the
representation theory of $\ih(n)$ can be found in \cite{MAT} and
some elements of the representation theory of $\mcalh(D_n)$ appear
for instance in \cite{Hu1} and in \cite{Hu2}. We call the
representation $\n^{(n)}$ introduced in \cite{CL} the Cohen-Wales
representation of the CGW algebra of type $D_n$. %We know from
%Theorem $2$, point $(ii)$ of \cite{CL} that for every $n\geq 4$,
%$\n^{(n)}$ is reducible if and only if $l\in\lb
%\unsur{r^{4n-7}},\unsur{r^{2n-7}},-\unsur{r^{2n-5}},r^3,\unsurr,-r^3\rb$.
We prove the following Theorem.

\newtheorem{theo}{Theorem}

\begin{theo}
Let $n$ be an integer with $n\geq 4$. Assume $\mcalh(D_n)$ and
$\ih(n)$ are both semisimple. Thus, assume that $r^{2k}\not\in\lb
-1,1\rb$ for every integer $k$ with $1\leq k\leq n-1$ and
$r^{2n}\neq 1$.\\
$1)$ \underline{Suppose first $n\geq 5$}. There are two cases.\\\\
(i) Assume $r^{2(n-1)}\not\in\lb i,-i\rb$ if $l=-r^3$. When the
Cohen-Wales representation $\n^{(n)}$ of the $CGW$ algebra of type
$D_n$ of degree $n^2-n$ is reducible, its unique proper invariant
subspace is isomorphic to one of the Specht modules
$$S^{(0),(n)},\;\; S^{(0),(n-1,1)},\;\; S^{(1),(n-1)},\;\; S^{(0),(n-2,2)},\;\; S^{(2),(n-2)},\;\; S^{(1),(n-2,1)},$$ which respectively arise if and
only if $$l=\unsur{r^{4n-7}},\;\; l=\unsur{r^{2n-7}},\;\;
l=-\unsur{r^{2n-5}},\;\; l=r^3,\;\;l=\unsurr,\;\; l=-r^3$$  $(ii)$
If $l=-r^3=\unsur{r^{4n-7}}$, there are exactly three proper
invariant subspaces in the Cohen-Wales space, and they are
respectively isomorphic to
$$S^{(0),(n)},\;S^{(1),(n-2,1)}\;\text{and}\;\;S^{(0),(n)}\oplus S^{(1),(n-2,1)}$$
$2)$\underline{Case $n=4$}. There are three cases.\\\\
(i) Assume that $r^6\neq i$ and $r^6\neq -i$ when $l=-r^3$. Assume
also that $l\neq\unsurr$. When the Cohen-Wales representation
$\n^{(4)}$ of the $CGW$ algebra of type $D_4$ of degree $12$ is
reducible, its unique proper invariant subspace is isomorphic to one
of the Specht modules
$$S^{(0),(4)},\;S^{(0),(2,2)},\;S^{(1),(3)},\;S^{(1),(2,1)},$$ which respectively arise if and only if
$$l=\unsur{r^9},\;l=r^3,\;l=-\unsur{r^3},\;l=-r^3$$
(ii) Same as (ii) in the general case $1)$.\\\\%If $l=r^3=\unsur{r^9}$,
%there are exactly three proper
%invariant subspaces in the Cohen-Wales space, and they are
%respectively isomorphic to
%$$S^{(0),(4)},\;S^{(1),(2,1)}\;\text{and}\;\;S^{(0),(4)}\oplus
%S^{(1),(2,1)}$$
(iii) if $l=\unsurr$, there are exactly seven proper invariant
subspaces in the Cohen-Wales space and they are respectively
isomorphic to
\begin{multline*}S^{(0),(3,1)},\;S^{(2,2)^{+}},\;S^{(2,2)^{-}},
S^{(0),(3,1)}\oplus\, S^{(2,2)^{+}},\;S^{(0),(3,1)}\oplus\,
S^{(2,2)^{-}},\;S^{(2,2)^{+}}\oplus\, S^{(2,2)^{-}},\\
S^{(0),(3,1)}\oplus\,
S^{(2,2)^{+}}\oplus\,S^{(2,2)^{-}}\end{multline*} A summary of the
Specht modules that occur in the Cohen-Wales space and the values of
the parameters for which they occur is given in Table $1$ below.
\end{theo}

\noindent Note $(ii)$ is the only case when two of the reducibility
values for $l$ can be equal in the case when $n\geq 5$. When $n=4$,
two of the reducibility values for $l$ in the generic case are
identical, which makes this case special and described in
$(iii)$.\\\\
\begin{tabular}[c]{|c|c|c|c|}
\hline$\begin{array}{l}\end{array}$&$\begin{array}{l}\end{array}$&$\begin{array}{l}\end{array}$&$\begin{array}{l}\end{array}$\\\text{Reducibility
value}&\text{Specht
module}&\text{Double partition}&\text{Dimension}\\$\begin{array}{l}\end{array}$&$\begin{array}{l}\end{array}$&$\begin{array}{l}\end{array}$&$\begin{array}{l}\end{array}$\\\hline$\begin{array}{l}\end{array}$&$\begin{array}{l}\end{array}$&$\begin{array}{l}\end{array}$&$\begin{array}{l}\end{array}$\\
$l=\unsur{r^{4n-7}}$&
$S^{(0),(n)}$&$\O\;,\;\;$\begin{tabular}[c]{|c|c|c|c|c|}\hline & & &
$\dots$ &
 \\\hline
\end{tabular}&$1$\\$\begin{array}{l}\end{array}$&$\begin{array}{l}\end{array}$&$\begin{array}{l}\end{array}$&$\begin{array}{l}\end{array}$\\\hline$\begin{array}{l}\end{array}$&$\begin{array}{l}\end{array}$&$\begin{array}{l}\end{array}$&$\begin{array}{l}\end{array}$\\
$l=\unsur{r^{2n-7}}$&$S^{(0),(n-1,1)}$&$\O\;,\;\;$\begin{tabular}[c]{|c|c|c|c|}\hline & &$\dots$&\\\hline\\\cline{1-1}\end{tabular}&$n-1$\\$\begin{array}{l}\end{array}$&$\begin{array}{l}\end{array}$&$\begin{array}{l}\end{array}$&$\begin{array}{l}\end{array}$\\\hline$\begin{array}{l}\end{array}$&$\begin{array}{l}\end{array}$&$\begin{array}{l}\end{array}$&$\begin{array}{l}\end{array}$\\
$l=-\unsur{r^{2n-5}}$&$S^{(1),(n-1)}$&\begin{tabular}[c]{|c|}\hline \\\hline\end{tabular}$\;,\;\;$\begin{tabular}[c]{|c|c|c|}\hline &$\dots$&\\\hline\end{tabular}&$n$\\$\begin{array}{l}\end{array}$&$\begin{array}{l}\end{array}$&$\begin{array}{l}\end{array}$&$\begin{array}{l}\end{array}$\\\hline$\begin{array}{l}\end{array}$&$\begin{array}{l}\end{array}$&$\begin{array}{l}\end{array}$&$\begin{array}{l}\end{array}$\\
$l=r^3$&$S^{(0),(n-2,2)}$&$\O\;,\;\;$\begin{tabular}[c]{|c|c|c|c|}\hline &&$\dots$&\\\hline &\\\cline{1-2}\end{tabular}&$\cil$\\$\begin{array}{l}\end{array}$&$\begin{array}{l}\end{array}$&$\begin{array}{l}\end{array}$&$\begin{array}{l}\end{array}$\\\hline$\begin{array}{l}\end{array}$&$\begin{array}{l}\end{array}$&$\begin{array}{l}\end{array}$&$\begin{array}{l}\end{array}$\\
$l=\unsurr$&$S^{(2),(n-2)}$&\begin{tabular}[c]{|c|c|}\hline&\\\hline\end{tabular}$\;,\;\;$\begin{tabular}[c]{|c|c|c|}\hline&\dots&\\\hline\end{tabular}&$\frac{n(n-1)}{2}$\\$\begin{array}{l}\end{array}$&$\begin{array}{l}\end{array}$&$\begin{array}{l}\end{array}$&$\begin{array}{l}\end{array}$\\\hline$\begin{array}{l}\end{array}$&$\begin{array}{l}\end{array}$&$\begin{array}{l}\end{array}$&$\begin{array}{l}\end{array}$\\
$l=-r^3$&$S^{(1),(n-2,1)}$&\begin{tabular}[c]{|c|}\hline
\\\hline\end{tabular}$\;,\;\;$\begin{tabular}[c]{|c|c|c|}\hline &\dots &\\\hline\\\cline{1-1}\end{tabular}&$n(n-2)$
\\$\begin{array}{l}\end{array}$&$\begin{array}{l}\end{array}$&$\begin{array}{l}\end{array}$&$\begin{array}{l}\end{array}$\\
\hline\end{tabular}\\ $\begin{array}{l}\\\end{array}$\\\textit{Table
$1$.}$\;$\textit{Classification of the invariant subspaces of the
Cohen-Wales representation of the CGW algebra of
type $D_n$}.\\

\noindent Our proof is based on the existing theorems of \cite{CL}.
We do the three following things. We identify the dimensions of the
irreducible $\mcalh(D_n)$-modules that may occur in the Cohen-Wales
space when $n\geq 11$, under some assumption at rank $10$. We prove
the classification theorem for these $n$, assuming the
classification theorem holds when $4\leq n\leq 10$ and we use in
particular some results of \cite{LEV}. We deal with the small cases
$4\leq n\leq 10$.
%holds in the small cases $n\in\lb 4,5,6,7,8,9,10\rb$. This allows us
%to rule out some dimensions for the invariant subspaces in the
%generic case.
\section{Proof of the classification theorem}
In \cite{CL}, $\S\,3.2$, we describe the degrees of the irreps of
$\mcalh(D_n)$ that are less than $n^2-n$, the degree of the
Cohen-Wales representation of type $D_n$. Our study combined with
Theorem $7$ of \cite{JAMES} shows that for $n\geq 9$, an irreducible
$\mcalh(D_n)$-module of dimension less than $n^2-n$ is one of the
Specht modules $S^{(0),(n)}$, $S^{(0),(n-1,1)}$, $S^{(0),(n-2,2)}$,
$S^{(0),(n-2,1,1)}$, $S^{(1),(n-1)}$, $S^{(1),(n-2,1)}$,
$S^{(2),(n-2)}$, $S^{(1,1),(n-2)}$, or one of their conjugates, or
is a Specht module $S^{(0),\la}$ for some partition $\la$ of $n$,
whose dimension lies between $\dbw$ and $n^2-n$. The Specht modules
that were listed first are of respective dimensions $1$, $n-1$,
$\cil$, $\dbw$, $n$, $n(n-2)$ and $\chl$ for the last two. Our goal
next is to show that the Specht modules of the second category
cannot occur in the Cohen-Wales space $V_n$ when $n\geq 11$, under
some assumption at rank $10$. We introduce some convenient
notations.
\newtheorem{Definition}{Definition}
\newtheorem{Proposition}{Proposition}
\newtheorem{Result}{Result}
\newtheorem{Lemma}{Lemma}
\begin{Definition}
We will denote by $Q_n(m)$ the set of Specht modules $S^{(0),\la}$
with the first row or first column of the Ferrers diagram of the
partition $\la$ of $n$ containing at least $n-m$ boxes.
\end{Definition}
\noindent When $n\geq 11$, we have seen at the end of $\S\,3.2$ of
\cite{CL} that there are no elements of $Q_n(3)\setminus Q_n(2)$
occuring in $V_n$ as their dimensions are too big. This is the key
remark to further show the following proposition under hypothesis
$\mathbf{(\mcalh)_{10}}$.\\\\
\textbf{$\mathbf{(\mcalh)_k}:$ The Specht modules $S^{(0),\la}$
occurring in $V_{k}$ belong to $Q_{k}(2)$.}
%We later prove that $(\mcalh)$ indeed holds.
\begin{Proposition}
%Assume $\mathbf{(\mcalh)_{10}}$ holds. Then, for any $n\geq 11$, the
%Specht modules $S^{(0),\la}$ occuring in $V_n$ belong to $Q_n(2)$.
$\mathbf{(\mcalh)_{10}}\Rightarrow \forall n\geq
11,\,\mathbf{(\mcalh)_{n}}$
\end{Proposition}
\noindent\textsc{Proof.} Let $n\geq 11$ and suppose the property
holds at rank $n-1$. Let $\W$ be a proper invariant subspace of
$V_n$ such that $\W$ is isomorphic to $S^{(0),\la}$ for some
partition $\la$ of $n$. Because $n\geq 11$, we already know that
$\W\not\in Q_n(3)\setminus Q_n(2)$. Suppose $\W\in Q_n(m)\setminus
Q_n(3)$ for some integer $m\geq 4$. Since in particular $\W\not\in
Q_n(2)$ and $n\geq 9$, Theorem $7$ of \cite{JAMES} implies that
$\text{dim}(\W)>\dbw$. And so $\W\cap V_{n-1}\neq\lb 0\rb$.
Moreover, we see with the branching rule that
$\W\da_{\mcalh(D_{n-1})}$ is isomorphic to a direct sum of an
element of $Q_{n-1}(m-1)\setminus Q_{n-1}(2)$ and an element of
$Q_{n-1}(m)\setminus Q_{n-1}(3)$. But by induction $\W\cap V_{n-1}$
belongs to $Q_{n-1}(2)$, hence a contradiction and the result.\\

\noindent The discussion at the beginning of $\S\,2$ and Proposition
$1$ imply Proposition $2$.\\
\begin{Proposition} Assume $(\mcalh)_{10}$ holds. Then, for any $n\geq 11$, the
Specht modules that may occur in the Cohen-Wales space $V_n$ have
dimensions
$$1,n-1,n,\frac{n(n-3)}{2},\frac{(n-1)(n-2)}{2},\frac{n(n-1)}{2},n(n-2)$$
\end{Proposition}
\noindent The discussion at the beginning of $\S\,2$, Proposition
$2$ and the theorems of the introduction of \cite{CL} and their
proofs imply in turn Proposition $3$.
\begin{Proposition}
Assume Theorem $1$ holds for integers $n$ with $4\leq n\leq 10$.
Then, for any $n\geq 4$, the Specht modules that may occur in the
Cohen-Wales space $V_n$ are
\begin{multline*}S^{(0),(n)}, S^{(0),(n-1,1)}, S^{(0),(n-2,2)},
S^{(0),(n-2,1,1)}, S^{(0),(3,1^{n-3})},
\\S^{(1),(n-1)}, S^{(1),(n-2,1)},
S^{(2),(n-2)}\end{multline*}
\end{Proposition}
\noin\textsc{Proof.} It suffices to show that the Specht modules
$S^{(1,1),(n-2)}$, $S^{(1,1),(1^{n-2})}$, $S^{(2),(1^{n-2})}$ when
$n\geq 5$, the Specht module $S^{(1,1),(2)}$ when $n=4$, and the
Specht module $S^{(1),(2,1^{n-3})}$
when $n\geq 5$ cannot occur in $V_n$.\\

First, we deal with the case $n=4$. We show the following result. We
use the notations of \cite{CL}, where the vectors $w_{ij}$'s and
$\wh{w_{ij}}$'s for $1\leq i<j\leq n$ denote the basis vectors of
the Cohen-Wales space. For the precise definition of the Cohen-Wales
representation of the CGW algebra of type $D_n$, we refer the reader
to Theorem $1$ of \cite{CL}.

\begin{Result}
$(i)$ There exists in $V_4$ an irreducible $3$-dimensional invariant
subspace isomorphic to $S^{(2,2)^{+}}$ (resp $S^{(2,2)^{-}}$) if and
only if $l=\unsurr$. If so, it is unique.\\
Moreover, the $6$-dimensional invariant subspace spanned by the
vectors $t_{ij}=w_{ij}-\wh{w_{ij}}$ with $1\leq i<j\leq 4$
decomposes as a direct sum of two irreducible invariant subspaces
$\U$ and $\W$, both of dimension $3$ and respectively spanned by the
vectors
$$\left\lb\begin{array}{ccc}
u_1&=&m(w_{12}-\wh{w_{12}})+(w_{13}-\wh{w_{13}})+(w_{24}-\wh{w_{24}})\\
u_2&=&(1-r^2)(w_{12}-\wh{w_{12}})+(w_{23}-\wh{w_{23}})+(w_{14}-\wh{w_{14}})-m\,(w_{24}-\wh{w_{24}})\\
u_3&=&(w_{12}-\wh{w_{12}})+(w_{34}-\wh{w_{34}})\end{array}\right.$$
$$\left\lb\begin{array}{ccc}
w_1&=&-(w_{13}-\wh{w_{13}})+m\,(w_{23}-\wh{w_{23}})+(w_{24}-\wh{w_{24}})\\
w_2&=&-(w_{12}-\wh{w_{12}})+m\,(w_{13}-\wh{w_{13}})+(1-r^2)(w_{23}-\wh{w_{23}})+(w_{34}-\wh{w_{34}})\\
w_3&=&-(w_{23}-\wh{w_{23}})+(w_{14}-\wh{w_{14}})\end{array}\right.$$
$(ii)$ The Specht modules $S^{(1^2,1^2)^{+}}$ and
$S^{(1^2,1^2)^{-}}$
don't occur in $V_4$. \\
$(iii)$ The Specht module $S^{(1,1),(2)}$ does not occur in $V_4$.
\end{Result}

\noin \textsc{Proof.} Assume $l=\unsurr$. Then, for any $n\geq 4$,
the $t_{ij}$'s span an $\frac{n(n-1)}{2}$-dimensional invariant
subspace of $V_n$. When $n\geq 5$, it is shown in $\S\,4$ of
\cite{CL} that this invariant subspace of $V_n$ is irreducible. We
show when $n=4$, the $6$-dimensional invariant subspace $\T$ spanned
by the $t_{ij}$'s is reducible and decomposes as a direct sum of two
irreducible $3$-dimensional invariant subspaces of $V_4$. We formed
the matrices $G_i$'s of the left actions by the $g_i$'s, $1\leq
i\leq 4$ in the basis $(w_{ij}-\wh{w_{ij}})_{1\leq i<j\leq 4}$ of
$\T$. Using Maple, we then computed the commutant of these matrices
and derived the following informations. The commutant has dimension
$2$ and is spanned by a homothety and another matrix $A$. The latter
matrix is symmetric and has two eigenvalues, namely
$$\la_1=\frac{1-r^2+2r^4}{r(1-r^2)}\qquad\;\text{
and}\;\qquad\la_2=\frac{r(r^4-r^2+2)}{r^2-1}$$ The corresponding
eigenspaces both have dimension $3$ and are respectively spanned by
the vectors $u_1$, $u_2$, $u_3$ and $w_1$, $w_2$, $w_3$ of $\T$
defined in the statement of the result. Since $A$ commutes to the
$G_i$'s, if $v_k$ is an eigenvector of $A$ for the eigenvalue
$\la_k$, then $G_iv_k$ is also an eigenvector of $A$ for the same
eigenvalue. So, the vector spaces $\U$ and $\W$ are invariant under
the actions by the $g_i$'s. We provide below the matrices $M_i$'s
(resp $N_i$'s) of the left actions by the $g_i$'s in the basis
$(u_1,u_2,u_3)$ (resp $(w_1,w_2,w_3)$) of $\U$ (resp $\W$).
$$\begin{array}{cc}M_2=\begin{bmatrix} 0
&1&0\\1&-m&0\\0&0&r\end{bmatrix}&\\&\\M3=\begin{bmatrix}
0&m\,r&1\\0&r&0\\1&-m&-m\end{bmatrix}&M1=\begin{bmatrix}
-m&-1&\frac{m}{r}\\-1&0&m\\0&0&r\end{bmatrix}\\&\\
M4=\begin{bmatrix}0&1&0\\1&-m&0\\0&0&r\end{bmatrix}&
\end{array}$$
\vspace{0.2cm} $\qquad\qquad\qquad$\textit{Matrix representation for
the} $\mcalh(D_4)$\textit{-module} $\U$.
$$\begin{array}{cc} N_2=\begin{bmatrix}
-m&0&1\\0&r&0\\1&0&0\end{bmatrix}&\\&\\N_3=\begin{bmatrix}0&1&0\\1&-m&0\\0&0&r\end{bmatrix}&
N_1=\begin{bmatrix}0&m\,r&-1\\0&r&0\\-1&m&-m\end{bmatrix}\\&\\
N_4=\begin{bmatrix}0&m\,r&-1\\0&r&0\\-1&m&-m\end{bmatrix}&\end{array}$$
\vspace{0.2cm} $\qquad\qquad\qquad$\textit{Matrix representation for
the} $\mcalh(D_4)$\textit{-module} $\W$.\\
\noin To show that $\U$ and $\W$ are irreducible, it suffices to
notice that the existence of a one-dimensional invariant subspace
inside them would force $l=\unsur{r^9}$ by Theorem $3$ of \cite{CL}.
This not compatible with $l=\unsurr$ and $r^8\neq 1$.\\
When $l=\unsurr$, there exists in $V_4$ a \textbf{unique}
irreducible $3$-dimensional invariant subspace $\V$ that is
isomorphic to $S^{(0),(3,1)}$ by the proof of Theorem $4$ of
\cite{CL} and there does not exist any irreducible $3$-dimensional
invariant subspace of $V_4$ that is isomorphic to $S^{(0),(2,1,1)}$.
Since $\U\neq\V$ and $\W\neq\V$, the two irreducible matrix
representations of $\mcalh(D_4)$ of degree $3$ defined by the
matrices $M_i$'s and $N_i$'s above are not equivalent to any of
those defined by the relations $(\Delta)$ and $(\nabla)$ on pages
$24$ and $25$ of \cite{CL}. Then they must be matrix representations
for $S^{(2,2)^{+}}$ or $S^{(2,2)^{-}}$. Moreover, using for instance
Maple, we see that the two matrix representations $(M_i)_{1\leq
i\leq 4}$ and $(N_i)_{1\leq i\leq 4}$ are inequivalent. Then, the
sufficient conditions
in point $(i)$ of Result $1$ hold.\\
To finish the proof of Result $1$, we prove the following
Proposition.
\begin{Proposition}\hfill\\
$(i)$ If there exists in $V_4$ an irreducible $3$-dimensional
invariant subspace, then $l=\unsurr$.\\
$(ii)$ If there exists in $V_4$ an irreducible $6$-dimensional
invariant subspace, then $l=\unsurr$.
\end{Proposition}
First we show a lemma.
\begin{Lemma}
If in the Cohen-Wales space $V_4$ there exists an irreducible
$2$-dimensional invariant subspace isomorphic to $S^{(0),(2,2)}$,
then $l=r^3$.
\end{Lemma}
\noindent A consequence of that lemma is that Theorem $6$ of
\cite{CL} holds in the case $n=4$ as well. To show the lemma, we
note that the Specht module $S^{(0),(2,2)}$ is the only irreducible
$\mcalh(D_4)$-module of dimension $2$ and we use the irreducible
matrix representation of degree $2$ of $\mcalh(D_4)$ that was given
on page $32$ of \cite{CL}. Suppose there exists two linearly
independent vectors $u_1$ and $u_2$ of
$V_4$ such that\\\\
$\forall \,i\in\lb 1,2,4\rb,\; \left\lb\begin{array}{ccc}g_i\,u_1&=&-\unsurr\,u_1\;\;\spadesuit_i\\
g_i\,u_2&=&u_1+r\,u_2\end{array}\right.\;$ and
$\;\left\lb\begin{array}{ccc}g_3\,u_1&=&r\,u_1+u_2\\g_3\,u_2&=&-\unsurr\,u_2\end{array}\right.$\\\\
Then, by $\spadesuit_2$ and $\spadesuit_4$,
\begin{equation}u_1=\ovl{\la_{13}}\,\ovl{w_{13}}-\unsurr\,\ovl{\la_{13}}\ovl{w_{14}}+\unsur{r^2}\,\ovl{\la_{13}}\ovl{w_{24}}
-\unsurr\,\ovl{\la_{13}}\ovl{w_{23}}\end{equation} Moreover, by
$\spadesuit_1$, we have $\la_{13}=\unsur{r^2}\,\wh{\la_{13}}$. In
particular, all the coefficients in $(1)$ are non-zero. Without loss
of generality, take $\la_{13}=1$. So, a complete expression for
$u_1$ is
\begin{equation}
u_1=w_{13}-\unsurr\,w_{14}+\unsur{r^2}\,w_{24}-\unsurr\,w_{23}+r^2\big(\wh{w_{13}}-\unsurr\,\wh{w_{14}}+\unsur{r^2}\,\wh{w_{24}}
-\unsurr\,\wh{w_{23}}\big)
\end{equation}
Now, the coefficient in $\wh{w_{12}}$ in $u_2=g_3.u_1-r\,u_1$ is
$r^2$. By looking at the coefficient of $\wh{w_{12}}$ in
$g_1\,u_2=u_1+r\,u_2$, we hence derive
$$\frac{r^4}{l}-\frac{m}{r^2}+m^2\bigg(\unsurr+r\bigg)=r^3$$
After simplification, this yields $l=r^3$.\\
In \cite{CL}, we introduce a $\cgwn$-submodule $K(n)$ of $V_n$ which
has the nice property that it contains all the proper invariant
subspaces of $V_n$ when the representation $\n^{(n)}$ is reducible.
We denoted its dimension by $k(n)$ (see Definition $3$ of
$\S\,3.6.1$ of \cite{CL}). In \cite{CL}, we use a program in
Mathematica that computes a matrix whose kernel is $K(n)$. Our
program runs until rank $7$ and computes $k(n)$ for the different
values of $l$ and $r$. By Theorem $2$ point $(i)$ of \cite{CL}, we
know that if $\n^{(4)}$ is reducible, then
$$l\in\bigg\lb \unsur{r^9},-r^3,r^3,-\unsur{r^3},\unsurr\bigg\rb$$
The corresponding values for $k(4)$ are in the same order
$$k(4)\in\lb 1,8,2,4,9\rb,$$
except when $l=\unsur{r^9}=-r^3$ when $k(4)=9$. Suppose there exists
in $V_4$ an irreducible $3$-dimensional invariant subspace. Then
$k(4)\geq 3$, and so $l\neq r^3$. Also, we have $l\neq
-\unsur{r^3}$. Indeed, otherwise the irreducible $3$-dimensional
invariant subspace would have a one-dimensional
$\mcalh(D_4)$-summand in $K(4)$. This would force $l=\unsur{r^9}$ by
Theorem $3$ of \cite{CL}. But we have $r^6\neq -1$ with our
restrictions on $r$. Suppose $l=-r^3$ and $r^{12}\neq -1$. Then we
have $k(4)=8$. Then, the irreducible $3$-dimensional invariant
subspace has a summand in $K(4)$ of dimension $5$. By the
representation theory of $\mcalh(D_4)$, this summand is not
irreducible. By Lemma $1$, this summand cannot contain an
irreducible $2$-dimensional invariant subspace. Nor can it contain
an irreducible $4$-dimensional invariant subspace by Theorem $5$,
point $(i)$ of \cite{CL}. So, it is impossible to have $k(4)=8$. We
cannot have $l=-r^3=\unsur{r^9}$ either. Indeed, if $k(4)=9$, the
irreducible $3$-dimensional invariant subspace has a summand in
$K(4)$ of dimension $6$. Moreover, since $l=\unsur{r^9}$, this
summand must contain a one-dimensional invariant subspace. Then
there exists in $V_4$ a $5$-dimensional invariant subspace and we
conclude as before, using also the uniqueness part in Theorem $3$ of
\cite{CL}.
The only possibility that is left is hence to have $l=\unsurr$. This finishes the proof of $(i)$. \\
Suppose there exists in $V_4$ an irreducible $6$-dimensional
invariant subspace. Name it $\V_6$. Then $k(4)\geq 6$ and so
$k(4)\in\lb 8,9\rb$. If $k(4)=9$ and $l=\unsurr$, then we are done.
Suppose $k(4)=9$ and $l=\unsur{r^9}=-r^3$. Then $\V_6$ has a
$3$-dimensional summand in $K(4)$. Since $l\neq r^3$, this summand
must be irreducible. Then an application of point $(i)$ yields
$l=\unsurr$, a contradiction. Finally, if $k(4)=8$, then $\V_6$ has
a $2$-dimensional summand in $K(4)$. By uniqueness of a
one-dimensional invariant subspace of $V_4$ when it exists, this
summand must be irreducible. It then follows from Lemma
$1$ that $l=r^3$. This is impossible.\\

\noin Using Proposition $4$, we close the proof of Result $1$ as
follows. For $(ii)$, we use the fact that when $\n^{(4)}$ is
reducible, it is indecomposable (this is a consequence of
Proposition $1$ of
$\S\,3.1$ of \cite{CL}). The same argument yields the uniqueness part in point $(i)$.\\
For $(iii)$, we use the fact that if $S^{(1,1),(2)}$ occurs in
$V_4$, then by Proposition $4$ point $(ii)$, we get $l=\unsurr$.
Then, the sum of the degrees of the irreducible
$\mcalh(D_4)$-modules occurring in $V_4$ exceeds the degree of the
Cohen-Wales representation of the CGW algebra of type $D_4$.
\\

\noin We state below a corollary of Result $1$ and of the proof of
Theorem $4$ of \cite{CL}.

\newtheorem{Corollary}{Corollary}
\begin{Corollary}
There exists in $V_4$ an irreducible $3$-dimensional invariant
subspace if and only if $l=\unsurr$. If so, $V_4$ contains exactly
three irreducible invariant subspaces of dimension $3$. They are
respectively spanned by the $u_i$'s and the $w_i$'s of Result $1$
and the $v_i$'s that were defined in Theorem $4$ of \cite{CL}.
\end{Corollary}

%\textsc{Proof.} Follows immediately from points $(i)$ and $(ii)$ of
%Result $1$ and from the proof of Theorem $4$ of \cite{CL}.
\noin We now continue the proof of Proposition $3$ by showing the
following result.

\begin{Result}
For any $n\geq 5$, the Specht modules $S^{(2),(1^{n-2})}$,
$S^{(1,1),(1^{n-2})}$ and $S^{(1),(2,1^{n-3})}$ cannot occur in
$V_n$.
\end{Result}

\noindent Indeed, let $n\geq 5$ and suppose $\W\subset V_n$ and
$\W\simeq S^{(2),(1^{n-2})}$. Then $$\W\da_{\mcalh(D_{n-1})}\simeq
S^{(1),(1^{n-2})}\oplus S^{(2),(1^{n-3})}$$ But by the key arguments
at the end of $\S\,3.5$ of \cite{CL}, the Specht module
$S^{(1),(1^{n-2})}$ cannot be a constituent of
$\W\da_{\mcalh(D_{n-1})}$. \\The argument
is identical for $S^{(1,1),(1^{n-2})}$ and for $S^{(1),(2,1^{n-3})}$, hence the result.\\

\noindent We further show the following Theorem, which gives a
necessary condition on $l$ and $r$ for the existence of an
irreducible $\chl$-dimensional invariant subspace of $V_n$ when
$n\geq 5$, under some conditions when $n=7$ or $n\geq 9$.
\vspace{0.2cm}
\begin{theo}
Let $n$ be an integer with $n\geq 5$. When $n=7$ or $n\geq 9$,
assume that $(\mcalh)_n$ holds. If in the Cohen-Wales space $V_n$
there exists an irreducible $\chl$-dimensional invariant subspace,
then $l=\unsurr$. Moreover, it is isomorphic to the Specht module
$S^{(2),(n-2)}$.
\end{theo}
\vspace{0.15cm} \noin\textsc{Proof.} The hypothesis ensures that
there are no irreducible Specht modules $S^{(0),\la}$ with $\la$ a
partition of $n$ of dimension $\chl$ occurring in $V_n$ when $n\geq
9$. As for $n=7$, it will become clear later why we don't need to
study the occurrence inside $V_7$ of $S^{(0),(3,3,1)}$ or its
conjugate Specht module, both of dimension $21$.

First, we deal with some special cases. We show that the Specht
module $S^{(0),(5,3)}$ of dimension $28$ cannot occur in $V_8$.
%Notice $28>26$, so that
If $\W$ is an invariant subspace of $V_8$ that is isomorphic to
$S^{(0),(5,3)}$, then
$$\W\da_{\mcalh(D_5)}\simeq 3\,S^{(0),(3,2)}\oplus
3\,S^{(0),(4,1)}\oplus S^{(0),(5)}$$ But the proof of Lemma $3$
point $(i)$ of \cite{CL} shows that if there exists an invariant
subspace of $\W\da_{\mcalh(D_5)}$ that is isomorphic to
$S^{(0),(3,2)}$, then it is unique. Hence a contradiction.

%$\W\cap V_6\neq\lb 0\rb$. Since $\W\cap V_7$
%cannot be isomorphic to $S^{(0),(4,3)}$ by Proposition $3$ of
%\cite{CL}, by the branching rule, this intersection must be
%isomorphic to $S^{(0),(5,2)}$. Then, by Lemma $3$ point $(i)$ of
%\cite{CL} applied with $n=7$, we get $l=r^3$. Further, the
%$\mcalh(D_6)$-module $\W\cap V_6$ cannot be isomorphic to
%$S^{(0),(5,1)}$, as otherwise we would have $l=\unsur{r^5}$ by
%Theorem $4$ of \cite{CL}. So instead, $\W\cap V_6\simeq
%S^{(0),(4,2)}$. It follows
%$$\W\cap V_6\da_{\mcalh(D_5)}\simeq S^{(0),(4,1)}\oplus
%S^{(0),(3,2)}$$ %Then, by the proof of Lemma $3$ point $(i)$ of
%\cite{CL}, there exists in $\W\cap V_6$ some vectors $v_1$, $v_2$,
%$v_3$, $v_4$ and $v_5$ as defined by equations $(37)-(41)$ of the
%same lemma. Since $\W$ is an invariant subspace of $V_8$, we have
%$(g_8g_7g_6g_5).v_5\in\W$. This implies
%$$r(w_{18}+r^2\,\wh{w_{18}})-(w_{28}+r^2\,\wh{w_{28}})+r^5(w_{23}+r^2\wh{w_{23}})-r^6(w_{13}+r^2\wh{w_{13}})\in
%K(8)$$ Recall from \cite{CL} that by definition, an element of
%$K(8)$ is annihilated by all the algebra elements $C_{ij}$'s and
%$\wh{C_{ij}}$'s with $1\leq i<j\leq 8$. These elements appear in
%Definition $1$ of $\S\,3.6.1$ of \cite{CL}. The actions of these
%elements on some of the basis vectors of the Cohen-Wales space are
%gathered in the tables of Proposition $2$ of $\S\,3.6.2$ of
%\cite{CL}.
%The rest of the proof is now identical to the one of Proposition $3$
%in $\S\,4$ of \cite{CL}. \\
We now eliminate in turn the conjugate Specht module
$S^{(0),(2^3,1^2)}$. If $\W$ is an invariant subspace of $V_8$
isomorphic to $S^{(0),(2^3,1^2)}$, we have
$$\W\da_{\mcalh(D_6)}\simeq 2\,S^{(0),(2^2,1^2)}\oplus
S^{(0),(2^3)}\oplus S^{(0),(2,1^4)}$$ Moreover, since $28>26$, we
have $\W\cap V_6\neq\lb 0\rb$. By Lemma $3$ point $(ii)$ of
\cite{CL} (resp the end of $\S\,3.4$ in \cite{CL} , resp the proof
of Theorem $4$ of \cite{CL}), the Specht module $S^{(0),(2^2,1^2)}$
(resp $S^{(0),(2^3)}$, resp $S^{(0),(2,1^4)}$) cannot occur in
$V_6$.
Hence we get a contradiction.\\
We deal with the last special case, namely $S^{(1),(2,2)}$. If $\W$
is an invariant subspace of $V_5$ that is isomorphic to
$S^{(1),(2,2)}$, then $$\W\da_{\mcalh(D_4)}\simeq
S^{(0),(2,2)}\oplus S^{(1),(2,1)}.$$ The presence of a fifth node
does not modify the proof of Lemma $1$, and we derive $l=r^3$. Then,
by using Mathematica, we have $k(5)=5$. But the existence of an
invariant subspace of $V_5$ of dimension $10$
forces $k(5)\geq 10$, hence a contradiction.\\

We now deal with the general case. We show that $S^{(1,1),(n-2)}$
cannot occur in $V_n$. \\%First, when $n=4$, if there
%exists in $V_4$ an irreducible $6$-dimensional invariant subspace,
%then $k(4)\geq 6$ and it must be isomorphic to $S^{(1,1),(2)}$.
%Mathematica gives the values of $k(4)$ for the different
%reducibility values for $l$ and $r$. Either $k(4)=8$ and $l=-r^3$ or
%$k(4)=9$ and $l=\unsurr$. Now, having $k(4)=8$ is impossible, as
%there would exist an irreducible $2$-dimensional invariant subspace
%in $V_4$. Then, by Proposition $5$ of \cite{CL}, we would have
%$l=r^3$. This contradicts $l=-r^3$. Hence, $k(4)=9$ and $l=\unsurr$.
Notice when $n\geq 5$, we have $\chl>2(n-1)$. So, if $\W\simeq
S^{(1,1),(n-2)}$, then $\W\cap V_{n-1}\neq\lb 0\rb$. First, we deal
with the cases $n=5$ and $n=6$, then we proceed by induction on
$n\geq 6$. Let $\W\subset V_5$ with $\W\simeq S^{(1,1),(3)}$. Then,
$\W\da_{\mcalh(D_4)}\simeq S^{(1),(3)}\oplus S^{(1,1),(2)}$. On the
one hand, by point $(iii)$ of Result $1$, we know that $\W\cap V_4$
cannot be isomorphic to $S^{(1,1),(2)}$. On the other hand, if
$\W\cap V_4\simeq S^{(1),(3)}$, then by Result $1$ of $\S\,3.5$ of
\cite{CL}, we get $l=-\unsur{r^3}$. As both $\n^{(4)}$ and
$\n^{(5)}$ are reducible, under our restrictions on $r$, we must
then have $l=-\unsur{r^3}=\unsur{r^{13}}$. By using Mathematica as
explained before, we computed $k(5)=1$ when $l=\unsur{r^{13}}$ and
$l\neq -r^3$. This
is not compatible with $\text{dim}(\W)=10$ and $\W\subseteq K(5)$. \\
Suppose now $\W$ is an invariant subspace of $V_6$ with $\W\simeq
S^{(1,1),(4)}$. And so $\W\da_{\mcalh(D_5)}\simeq
S^{(1,1),(3)}\oplus S^{(1),(4)}$. Suppose $\W\cap V_5\simeq
S^{(1),(4)}$. Then, by using the same results as before, we get
$l=-\unsur{r^5}=\unsur{r^{17}}$. Again, the contradiction comes from
$k(6)=1$ when $l=\unsur{r^{17}}$ and $l\neq -r^3$. Hence $\W\cap
V_5$ is isomorphic to $S^{(1,1),(3)}$ instead. By the case
$n=5$, this is impossible.\\
Let $n\geq 7$ and suppose $S^{(1,1),(n-3)}$ does not occur in
$V_{n-1}$. Let $\W\subset V_n$ with $\W\simeq S^{(1,1),(n-2)}$.
Then,
$$\W\da_{\mcalh(D_{n-1})}\simeq S^{(1),(n-2)}\oplus S^{(1,1),(n-3)}$$
Since
$$\text{dim}(\W\cap V_{n-1})\geq \chl-2(n-1)=\frac{(n-1)(n-4)}{2},$$
we see that $\text{dim}(\W\cap V_{n-1})>n-1$ as soon as $n\geq 7$.
This contradicts the fact that $S^{(1,1),(n-3)}$ does not occur in
$V_{n-1}$. \\
It remains to show that if $\W$ is an invariant subspace of $V_n$
that is isomorphic to $S^{(2),(n-2)}$, then $l=\unsurr$. Again we
proceed by induction on $n\geq 5$ and adapt the proof above. Suppose
$\W$ is an invariant subspace of $V_5$ that is isomorphic to
$S^{(2),(3)}$. By the same arguments as above, the Specht module
$S^{(1),(3)}$ cannot occur in $\W\cap V_4$. Hence $S^{(2,2)^{+}}$ or
$S^{(2,2)^{-}}$ or both Specht modules must occur in $\W\cap V_4$.
Then, by Result $1$ point $(i)$, it follows that $l=\unsurr$.
Suppose now $\W$ is an invariant subspace of $V_6$ that is
isomorphic to $S^{(2),(4)}$. Then $\W\cap V_5$ must be isomorphic to
$S^{(2),(3)}$. Then, by the case $n=5$, we get $l=\unsurr$. \\
Finally, when $n\geq 7$, if $\W$ is an invariant subspace of $V_n$
that is isomorphic to $S^{(2),(n-2)}$, by the same arguments as
above, $S^{(2),(n-3)}$ occurs in $V_{n-1}$ and by induction this
forces $l=\unsurr$. This finishes the proof of the theorem. \\

We now better Proposition $2$ and Proposition $3$ by showing Result
$3$ and Result $4$
\begin{Result}
Let $n\geq 4$. The Specht module $S^{(0),(n-2,1,1)}$ cannot occur in
the Cohen-Wales space $V_n$.
\end{Result}

\textsc{Proof.} When $n=4$, this result is part of the proof of
Theorem $4$ of \cite{CL} in the case $n=4$. Assume now $n\geq 5$. We
will use the results of \cite{LEV} to give a matrix representation
for $S^{(0),(n-2,1,1)}$. In \cite{LEV} $\S\,2$, we give a new
representation of the braid group on $n$ strands that is equivalent
to the Lawrence--Krammer representation. In that paper, the vectors
$u_1$, $u_2$, $u_3$ of Theorem $3.5$ and the vectors $V_s^{(k)}$
with $5\leq k\leq n,\,1\leq s\leq k-2$ of Theorem $3.15$ form a
basis $\B$ of vectors of the Lawrence--Krammer space. It is a result
of the paper that the matrices of the left braid group actions by
the $g_i$'s, $1\leq i\leq n-1$, in this basis provide a matrix
representation for the $\ih(n)$-Specht module $S^{(n-2,1,1)}$. The
left braid group actions by the $g_i$'s on the basis of the
Lawrence--Krammer space is the braid group representation provided
in $\S\,2$ of \cite{LEV}. The braid group actions by the $g_i$'s on
the family of vectors $(V_s^{(k)})_{5\leq k\leq n,\, 1\leq s\leq
n-2}$ is the object of Lemma $3.16$ of \cite{LEV}. Moreover, we show
along the proof of Theorem $3.15$ in \cite{LEV} that if $r^8\neq
-1$, another equivalent irreducible matrix representation for
$S^{(n-2,1,1)}$ is provided by the matrices of the left braid group
actions by the $g_i$'s, $1\leq i\leq n-1$ in the basis $\B^{'}$
consisting this time of all the vectors $V_i^{(k)}$ with $3\leq
k\leq n$ and $1\leq i\leq k-2$. Denote by the $G_i$'s the matrices
of the left braid group actions by the $g_i$'s in $\B^{'}$. Define
$H_1=H_2=G_1$ and $H_i=G_{i-1}$ for every integer $i$ with $3\leq
i\leq n$. Since $n\geq 5$, our restrictions on $r$ impose in
particular $r^8\neq -1$. Then, the $H_i$'s define a matrix
representation for $S^{(0),(n-2,1,1)}$. If $S^{(0),(n-2,1,1)}$
occurs in the Cohen--Wales space, then there exists a basis of
$\frac{(n-1)(n-2)}{2}$ vectors
$$W_1^{(3)},W_1^{(4)},W_2^{(4)},W_1^{(5)},W_2^{(5)},W_3^{(5)},\dots,W_1^{(n)},\dots,W_{n-2}^{(n)}$$
of the Cohen--Wales space such that the matrices of the left
Cohen--Wales actions by the $g_i$'s, $1\leq i\leq n$ in this basis
are the $H_i$'s, $1\leq i\leq n$. The left Cohen--Wales action by
the $g_i$'s, $1\leq i\leq n$ on the basis $(\ovl{w_{ij}})_{1\leq
i<j\leq n}$ of the Cohen--Wales space is the Cohen-Wales
representation constructed and defined in Theorem $1$ of \cite{CL}
and whose structure is studied in this paper.

\begin{Lemma}
The vectors $W_i^{(k)}$'s, $3\leq k\leq n$, $1\leq i\leq k-2$, do
not contain any hat terms.
\end{Lemma}
\textsc{Proof.} First, we show that this is true for the vectors
$W_{j-2}^{(j)}$ with $3\leq j\leq n$. Second, we show that the
result also holds for all the other vectors $W_k^{(j)}$ with $3\leq
j\leq n$ and $1\leq k\leq j-3$. To do so, we proceed by descending
induction on the integer $k$.\\
We use the second and fourth equalities of Lemma $3.16$ of
\cite{LEV} to derive the relations
$$\left\lb\begin{array}{ccc}
g_{j-1}.\,W_{j-2}^{(j)}&=&-\unsurr\,W_{j-2}^{(j)}\\
&&\\ g_j.\,W_{j-2}^{(j)}&=&-\unsurr\,W_{j-2}^{(j)}
\end{array}\right.$$

\noindent These two relations imply that $W_{j-2}^{(j)}$ is a linear
combination of basis vectors from the Cohen-Wales space, such that
their nodes either start or end in $j-1$ or start in $j-2$ and end
in $j$. Further, the first (resp second) relation implies that there
is no term in $\wh{w_{j-2,j-1}}$ (resp $\wh{w_{j-1,j}}$) in
$W_{j-2}^{(j)}$. Furthermore, by the way the coefficients are
related as in Lemma $2$ of $\S\,3.3$ of \cite{CL} with
$\la=-\unsurr$, we see that $W_{j-2}^{(j)}$ is a multiple of
$$w_{j-2,j-1}-\unsurr\,w_{j-2,j}+\unsur{r^2}\,w_{j-1,j}$$
Fix $j\geq 4$ and suppose the vector $W_k^{(j)}$ does not contain
any hat term, where $2\leq k\leq j-2$. We show this implies that
$W_{k-1}^{(j)}$ does not contain any hat term either. From the first
equation in Lemma $3.16$ of \cite{LEV}, we derive
$$g_k.\,W_k^{(j)}=W_{k-1}^{(j)}+r\,W_k^{(j)}-r^{j-k-1}\,W_{k-1}^{(k+1)}\;\;\;\;\;\;(\star)$$
By the defining relations of the representation in Theorem $1$ of
\cite{CL}, if the vector $W_k^{(j)}$ does not contain any hat terms,
since $k\geq 2$, the vector $g_k.\,W_k^{(j)}$ does not contain any
hat term either. Moreover, by the first step, $W_{k-1}^{(k+1)}$ does
not contain any hat term. Then, with $(\star)$, we conclude that
$W_{k-1}^{(j)}$ does not contain any hat term. By induction, the
statement of the lemma holds for all the $W_s^{(j)}$'s with $1\leq
s\leq j-2$.\\
Using Lemma $1$, it is now easy to conclude. Indeed, take for
instance without loss of generality
$$W_1^{(3)}=w_{12}-\unsurr\,w_{13}+\unsur{r^2}\,w_{23},$$
and notice an action to the left by $g_1$ on $W_1^{(3)}$ creates a
term in $\wh{w_{23}}$ with coefficient $-\unsurr$. Then, by the
lemma, $g_1.\,W_1^{(3)}$ cannot be a linear combination of
$W_s^{(k)}$'s, which constitutes a contradiction. We conclude that
$S^{(0),(n-2,1,1)}$ cannot occur in $V_n$.

\begin{Result}
Let $n\geq 5$. The Specht module $S^{(0),(3,1^{n-3})}$ cannot occur
in the Cohen-Wales space $V_n$.
\end{Result}

\textsc{Proof.} Suppose $\W$ is an invariant subspace of $V_n$ that
is isomorphic to the Specht module $S^{(0),(3,1^{n-3})}$. Notice
when $n\geq 7$, the dimension of $\W$ is greater than $2(n-1)$ and
so $\W\cap V_{n-1}\neq\lb 0\rb$. Moreover, by the beginning of
$\S\,3.4$ of \cite{CL}, we know that the Specht module
$S^{(0),(2,1^{n-3})}$ does not occur in $V_{n-1}$ for any $n\geq 5$.
So, whenever $n\geq 7$, we get $\W\cap V_{n-1}\simeq
S^{(0),(3,1^{n-4})}$. Thus, by induction, the cases $n\geq 7$ reduce
to the case $n=6$. Let $\W\subset V_6$ such that $\W\simeq
S^{(0),(3,1^3)}$. If $\W\cap V_5\neq \lb 0\rb$, then $\W\cap
V_5\simeq S^{(0),(3,1,1)}$ and so the case $n=6$ reduces to the case
$n=5$. This case is treated below. We show that the intersection
$\W\cap V_5$ is indeed nonzero. %Otherwise, if $\W\cap V_5=\lb
%0\rb$, observe $\text{dim}(\W)+\text{dim}(V_5)=\text{dim}(V_6)$.
%Then, $\W\oplus V_5=V_6$.
An irreducible matrix representation for $S^{(0),(3,1^3)}$ is the
one given by the matrices $H_i$'s introduced in the proof of Result
$3$, where $r$ has been replaced by $-\unsurr$ and where $n=6$.
Denote these new matrices by the $K_i$'s. So, there exists a basis
$(U_i^{(k)})_{1\leq k\leq 6,1\leq i\leq k-2}$ of vectors of $\W$
such that the matrices of the left Cohen-Wales actions by the
$g_i$'s in this basis are the $K_i$'s. From the last equation of
Lemma $3.16$ of \cite{LEV}, we derive
$$g_i.\,U_s^{(k)}=-\unsurr\,U_s^{(k)}\qquad\forall i\not\in\lb
s,s+1,s+2,k\rb$$ In particular, we have
$$g_i.\,U_4^{(6)}=-\unsurr\,U_4^{(6)}\qquad\forall i\leq 3$$
This implies $U_4^{(6)}$ is a linear combination of $w_{12}$,
$w_{23}$ and $w_{13}$. So, $U_4^{(6)}$ belongs to $\W\cap V_5$.

Thus, it suffices to deal with the case $n=5$. We must show that it
is impossible to have an invariant subspace of $V_5$ that is
isomorphic to the Specht module $S^{(0),(3,1,1)}$. Observe this is
Result $3$
for $n=5$.\\

\noindent With Result $3$ and Result $4$, Proposition $2$ and
Proposition $3$ simplify nicely and become Proposition $5$ and
Proposition $6$ respectively.
\begin{Proposition}
Assume $(\mcalh)_{10}$ holds. Then, for any $n\geq 11$, the Specht
modules that may occur in the Cohen-Wales space $V_n$ have
dimensions
$$1,n-1,n,\frac{n(n-3)}{2},\frac{n(n-1)}{2},n(n-2)$$
\end{Proposition}

\begin{Proposition}
Assume Theorem $1$ holds for integers $n$ with $4\leq n\leq 10$.
Then, for any $n\geq 4$, the Specht modules that may occur in the
Cohen-Wales space $V_n$ are $$S^{(0),(n)}, S^{(0),(n-1,1)},
S^{(0),(n-2,2)}, S^{(1),(n-1)}, S^{(1),(n-2,1)}, S^{(2),(n-2)}$$
\end{Proposition}

%In what follows, we assume Theorem $1$ holds for the integers $n$
%with $4\leq n\leq 10$ and we show Theorem $1$ for the integers $n$
%such that $n\geq 11$.

\noindent We now show even more properties on the representation. We
prove in particular that point $(i)$ of Theorem $5$ of \cite{CL} is
in fact an equivalence. We have the following result.
\begin{theo}
Let $n\geq 4$ and suppose $l=-\unsur{r^{2n-5}}$. Then, $K(n)$ is the
unique invariant subspace of $V_n$ and $k(n)=n$. In particular,
there exists in $V_n$ a unique irreducible $n$-dimensional invariant
subspace.
\end{theo}

\noindent \textsc{Proof.} With Mathematica, we check that for these
values of $l$ and $r$, we have $k(n)=n$ for every $4\leq n\leq 7$.
Suppose $n\geq 8$. Then, by using points $(i)$ and $(iii)$ of
Theorem $8$ in \cite{CL} on one hand, Theorems $3$ and $4$ of
\cite{CL} on the other hand and with our restrictions on $r$, we see
that $k(n)\geq n$. When $n\geq 7$, we have $2n\leq \cil$. Next, if
$k(n)\geq 2n$, then $k(n)>2n-2$ and so $K(n)\cap V_{n-1}\neq \lb
0\rb$. This is impossible, because for these values of $l$ and $r$,
the representation $\n^{(n-1)}$ is irreducible. When $n=8$, one
cannot have $k(8)>14$ by the same argument. Moreover, if $k(8)=14$,
then $K(8)$ is irreducible and $K(8)\simeq S^{(0),(4,4)}$ or
$K(8)\simeq S^{(0),(2^3)}$. In the first situation, we get $l=r^3$
by the proof of Lemma $3$, point $(i)$ of \cite{CL}. This
contradicts $l=-\unsur{r^{11}}$ and $r^{2(8-1)}\neq -1$. The second
situation is impossible by the proof of Lemma $3$, point $(ii)$ of
\cite{CL}. Thus, when $l=-\unsur{r^{2n-5}}$, we have $k(n)=n$ for
all $n$. Then, $K(n)$ is irreducible, hence is the unique invariant
subspace of $V_n$.\\

We now show that the statement of Theorem $6$ of \cite{CL} is in
fact an equivalence for $n=4$ and $n\geq 6$ under $\h{8}$, $\h{9}$
and $\h{10}$ (Recall Theorem $6$ of \cite{CL} also holds for $n=4$
as shown by Lemma $1$ of the current paper).

\begin{theo}
Let $n\geq 4$. Suppose $\h{8}$, $\h{9}$ and $\h{10}$ hold. If
$l=r^3$, then $k(n)=\cil$ and $K(n)$ is the unique invariant
subspace of $V_n$. In particular, there exists in $V_n$ a unique
irreducible $\cil$-dimensional invariant subspace.
\end{theo}

\textsc{Proof.} The equality on the dimension is true for $4\leq
n\leq 7$ by using Mathematica. Suppose $n\geq 8$ and suppose
$k(n-1)=\frac{(n-1)(n-4)}{2}$. By Theorems $3,4,5$ of \cite{CL} and
Theorem $8$, points $(i)$ and $(iii)$ of \cite{CL}, we have
$k(8)\geq 14$ and $k(n)\geq\cil$ when $n\geq 9$. Further, under the
assumptions of the theorem, by Proposition $1$, by the study in
$\S\,3.2$ of \cite{CL}, by Results $3$ and $4$ above and by Theorem
$2$ of this paper, we see that the irreducible invariant subspaces
that may occur in $V_n$ have dimensions $\cil$, or $n(n-2)$, except
when $n=8$ when they can also have dimension $35$. We first deal
with the case $n\geq 9$. Suppose $k(n)\geq n(n-3)$. Then we have
$$k(n-1)\geq\text{dim}(K(n)\cap V_{n-1})\geq n(n-3)-2(n-1)$$
Observe the member to the right of the second inequality is greater
than $\frac{n(n-4)}{2}$ as soon as $n\geq 6$. This contradicts
$k(n-1)=\frac{(n-1)(n-4)}{2}$. \\
Suppose now $n=8$. It is shown in forthcoming Result $5$ that there
does not exist in $V_8$ any irreducible invariant subspaces of
dimension $14$. Hence, like in the general case, we get $k(8)\geq
20$ and the dimensions of the irreducible invariant subspaces that
may occur in $V_8$ are the following: $20$, $35$, $48$. Like in the
general case, it is impossible to have $k(8)\geq 40$. Hence, suppose
$k(8)=35$. But then $$k(7)\geq\text{dim}(K(8)\cap V_7)\geq
35-2\times 7=21,$$
which contradicts $k(7)=14$. Thus, we have $k(8)=20$. \\
Hence, by induction, $k(n)=\cil$ for all $n\geq 4$. Consequently
also, the module $K(n)$ is irreducible. Therefore, $K(n)$ is the
unique invariant subspace of $V_n$ when $l=r^3$.\\

\noindent The next theorem studies reducibility in the case $l=-r^3$
with some restrictions.

\begin{theo}
Assume $(\mcalh)_{10}$ holds. Let $n$ be an integer with $n\geq 11$.
Suppose $l=-r^3$ and $r^{2(n-1)}\not\in\lb i,-i\rb$. Then, there
exists in $V_n$ an irreducible $n(n-2)$-dimensional invariant
subspace.
\end{theo}

\textsc{Proof.} Immediate by Proposition $5$ of the current paper,
Theorems $3$, $4$, $5$, $6$ of \cite{CL} and Theorem $2$ of the
current paper.

\begin{theo}
Let $n\geq 4$. Suppose $\h{8}$, $\h{9}$ and $\h{10}$ hold. If there
exists in $V_n$ an irreducible $n(n-2)$-dimensional invariant
subspace, then $l=-r^3$.
\end{theo}

\textsc{Proof.} When $n=4$, Mathematica gives the value of $k(4)$
depending on the values for $l$ and $r$. We deduce that if there
exists an irreducible $8$-dimensional invariant subspace in $V_4$,
the only possibility is to have $l=-r^3$. The case $n=5$ is also
done with Mathematica and is similar. Assume now $n\geq 6$. Since
when $n\geq 5$, we have $n(n-2)>2(2n-3)$, the existence of an
irreducible $n(n-2)$-dimensional invariant subspace of $V_n$ implies
that $\n^{(n-1)}$ and $\n^{(n-2)}$ are both reducible. Then, we must
have $l\in\lb r^3,-r^3,\unsurr\rb$. Further, by Theorem $7$ of
\cite{CL}, when $l=\unsurr$, there exists an irreducible
$\chl$-dimensional invariant subspace in $V_n$. But
$n(n-2)+\frac{n(n-1)}{2}>n(n-1)$, so this is impossible.
Furthermore, by Theorem $4$ above, when $l=r^3$, there exists in
$V_n$ an irreducible $\frac{n(n-3)}{2}$-dimensional invariant
subspace. This is where the assumptions in the statement of Theorem
$6$ become relevant. Observe when $n\geq 6$, we have
$\frac{n(n-3)}{2}>n$. Thus, we get a contradiction since $n(n-2)+n$
is the degree of the Cohen-Wales representation. We conclude that
$l=-r^3$.

\noindent We are now ready to conclude when $n\geq 11$, assuming
Theorem $1$ holds for the small values of $n$. We have the following
statement.

\begin{Proposition}
Assume Theorem $1$ holds for $4\leq n\leq 10$. Then Theorem $1$
holds for every $n\geq 11$.
\end{Proposition}
\noindent\textsc{Proof.} For point (i), immediate by gathering all
the results of \cite{CL} and of this paper.\\
For point $(ii)$, it remains to show that there exists in $V_n$ an
irreducible $n(n-2)$-dimensional invariant subspace. %When $n\in\lb
%4,5\rb$, the respective values for $k(4)$ and $k(5)$ obtained with
%Mathematica force the result. Assume now $n\geq 6$. Then,
By Theorem $10$ point $(ii)$ of \cite{CL}, we know that $K(n)\cap
V_{n-1}\neq \lb 0\rb$ since there is an element of $V_5$ that
belongs to $K(n)$ for all $n\geq 5$. Moreover, by point $(i)$ of
Theorem $1$, we know that $K(n-1)$ is irreducible and that
$k(n-1)=(n-1)(n-3)$. Then, $K(n)\cap V_{n-1}=K(n-1)$ and $K(n)$
cannot be one-dimensional. By Proposition $5$ and Theorem $2$ of the
current paper, and by Theorems $3,4,5,6$ of \cite{CL}, we must have
$k(n)=1+n(n-2)$. It yields the result.
%\newtheorem{Remark}{Remark}
%\begin{Remark}
%In Proposition $6$, it suffices to assume that $\h{8}$, $\h{9}$ and
%$\h{10}$ hold.
%\end{Remark}

\noindent It remains to show Theorem $1$ indeed holds when $4\leq
n\leq 10$. We will show that by excluding a few Specht modules
$S^{(0),\la}$, we can in fact exclude many more. This is the meaning
of Proposition $8$ below. First, we need to establish the following
result.

\begin{Result}
In $V_8$, the Specht modules $S^{(0),(4,4)}$ and $S^{(0),(2^4)}$
don't occur.
\end{Result}

\textsc{Proof.} The proof is identical to the one of Proposition $3$
in $\S\,4$ of \cite{CL}.

\begin{Proposition}
If $\h{6}$ holds and $Q_n(3)\setminus Q_n(2)=\emptyset$ when $7\leq
n\leq 10$, then $\h{k}$ holds for every $k\geq 4$.
\end{Proposition}

\textsc{Proof.} By the same arguments as in the proof of Proposition
$1$ , it suffices to show that if a module $\W$ belongs to
$Q_n(m)\setminus Q_n(3)$ with $7\leq n\leq 10$ and $m\geq 4$, then
$\W\cap V_{n-1}\neq\lb 0\rb$. Applying James'theory in \cite{JAMES},
when $n\geq 9$, an irreducible $\mcalh(D_n)$-module $S^{(0),\la}$
either belongs to $Q_n(2)$ or has dimension greater than $\dbw$.
When $n\geq 6$, we have $\dbw\geq 2(n-1)$, so this settles $n\in\lb
9,10\rb$. When $n=7$, the same statement holds with the exception of
$S^{(0),(4,3)}$ and its conjugate Specht module. But both modules
belong to $Q_7(3)$, which by assumption is forbidden. Finally, when
$n=8$, the statement above is not true because of the Specht module
$S^{(0),(4,4)}$ and its conjugate, both of dimension $14$. But we
have seen in Result $5$ that they cannot occur. So we are done with
all the cases.

\begin{Lemma}
If $S^{(0),(3,2,1)}$ does not occur in $V_6$, then $\h{6}$ holds.
\end{Lemma}
\textsc{Proof.} By the end of $\S\,3.4$ in \cite{CL}, the Specht
module $S^{(0),(3,3)}$ and its conjugate $S^{(0),(2^3)}$ cannot occur in $V_6$.\\

Using Proposition $8$ and Lemma $3$, we derive Proposition $9$.

\begin{Proposition} (Sketch of the final stage of the proof).
To finish the proof of Theorem $1$, it suffices to show that the
following list of Specht modules cannot occur in the Cohen-Wales
space $V_n$. We also provide their dimensions and recall the space
in which they may have lived. The integer $k$ denotes the size of
the first partition in the double partition.
$$\begin{array}{cc}
\hspace{-2.6cm}k=0&
\begin{array}{cccccccc}
S^{(0),(3,2,1)}&\begin{array}{l}S^{(0),(4,2,1)}\\S^{(0),(3,2,1^2)}\end{array}&\begin{array}{l}S^{(0),(6,3)}\\S^{(0),(2^3,1^3)}\end{array}
&\begin{array}{l}S^{(0),(7,3)}\\S^{(0),(2^3,1^4)}\end{array}&S^{(0),(4,1^3)}&\begin{array}{l}S^{(0),(5,1^3)}\\S^{(0),(4,1^4)}\end{array}
&\begin{array}{l}S^{(0),(6,1^3)}\\S^{(0),(4,1^5)}\end{array}&\begin{array}{l}S^{(0),(7,1^3)}\\S^{(0),(4,1^6)}\end{array}\\
16&35&48&75&20&35&56&84\\
V_6&V_7&V_9&V_{10}&V_7&V_8&V_9&V_{10}\end{array}\\\hspace{-2.6cm}&\\
\hspace{-2.6cm}k=1&
\begin{array}{l}S^{(1),(3,3)}\\S^{(1),(2^3)}\\\;\;\;\;\;35\\\;\;\;\;\;V_7\end{array}\\\hspace{-2.6cm}&\\
\hspace{-2.6cm}k=3&
\begin{array}{ccccc}S^{(3),(1^3)}&\begin{array}{l}S^{(3,3)^{+}}\\S^{(1^3,1^3)^{+}}\end{array}&\begin{array}{l}
S^{(3,3)^{-}}\\S^{(1^3,1^3)^{-}}\end{array}&\begin{array}{l}S^{(3),(4)}\\S^{(1^3),(1^4)}\end{array}
&\begin{array}{l}S^{(3),(1^4)}\\S^{(1^3),(4)}\end{array}\\
20&10&10&35&35\\
V_6&V_6&V_6&V_7&V_7\end{array}\\\hspace{-2.6cm}&\\
\hspace{-2.6cm}k=4&
\begin{array}{cc}\begin{array}{l}S^{(4,4)^{+}}\\S^{(1^4,1^4)^{+}}\end{array}&\begin{array}{l}S^{(4,4)^{-}}\\S^{(1^4,1^4)^{-}}\end{array}\\
35&35\\
V_8&V_8\end{array}\end{array}$$
\end{Proposition}

\noindent\textsc{Proof.} The dimensions of the modules of
$Q_n(3)\setminus Q_n(2)$ are provided in \cite{CL} at the end of
$\S\,3.2$. The dimension $\frac{n(n-2)(n-4)}{3}$ of
$S^{(0),(n-3,2,1)}$ and its conjugate Specht module is too large as
soon as $n\geq 8$. Next, we have $\frac{n(n-1)(n-5)}{6}\geq n(n-1)$
if and only if $n\geq 11$. So $S^{(0),(n-3,3)}$ and its conjugate
must be considered when $n\in\lb 9,10\rb$. When $n=7$, it is shown
in \cite{CL} $\S\,4$ that $S^{(0),(4,3)}$ and $S^{(0),(2^3,1)}$
cannot occur. This is Proposition $3$ of \cite{CL}. When $n=8$, we
saw along the proof of Theorem $2$ that $S^{(0),(5,3)}$ and
$S^{(0),(2^3,1^2)}$ cannot occur in $V_8$. As for
$S^{(0),(n-3,1^3)}$ and its conjugate, again
$\frac{(n-1)(n-2)(n-3)}{6}\geq n(n-1)$ if and only if $n\geq 11$, so
these Specht modules must be considered for all the values $7\leq
n\leq 10$. Once we have excluded all these modules from the row
$k=0$, hypothesis $\h{k}$ holds for all $k\geq 4$ by Proposition
$8$. In particular, $\h{8}$, $\h{9}$ and $\h{10}$ hold. So, Theorems
$2,4,6$ become true for every integer $n\geq 4$ without any
assumption. Also, point $(i)$ of Theorem $1$ now holds for every
integer $n\geq 11$ without any assumption. Further, if we succeed to
exclude the modules from the rows $k=1$, $k=3$ and $k=4$, which
appear to be the only potential candidates which have not yet been
studied by $\S\,3.2$ of \cite{CL}, then Proposition $5$ becomes true
for every integer $n\geq 4$. Then, we see that Theorem $5$ holds for
every integer $n\geq 4$ under the only assumption that
$r^{2(n-1)}\not\in\lb i,-i\rb$. Then, by gathering all the results,
the proof of Theorem $1$, point $(i)$ is complete. As for $(ii)$,
the proof is the same as in Proposition $7$ when $n\geq 6$. And when
$n= 4$ (resp $n=5$), the value of $k(4)$ (resp $k(5)$) obtained with
Mathematica forces the existence of an irreducible invariant
subspace of $V_4$ (resp $V_5$) of dimension $8$ (resp $15$). The
proof is then complete.

\begin{Proposition} (End of the proof). The Specht modules from the
list of Proposition $8$ cannot occur in the Cohen-Wales
representation.
\end{Proposition}

\textsc{Proof.} We first deal with the row $k=0$. %Recall that if
%$\W$ is a proper invariant subspace of $V_n$ such that
%$\text{dim}(\W)>2(2n-3)$, then the intersection $\W\cap V_{n-2}$ is
%non-zero. Hence $\n^{(n-2)}$ is reducible. So is $\n^{(n-1)}$. Then
%by confronting the values for $l$ and $r$ in Theorem $2$, point
%$(i)$ of \cite{CL}, we get $l\in\lb\unsurr, r^3,-r^3\rb$. The table
%below gives the value of $2(2n-3)$ for the different values of
%$n$.\\
%\begin{center}
%\begin{tabular}{|c|c|c|c|c|}
%\hline $n=6$&$n=7$&$n=8$&$n=9$&$n=10$\\\hline 18&22&26&30&34\\\hline
%\end{tabular}
%\end{center}
%\vspace{0.3cm}With this table, we see that all the Specht modules
%present in the row $k=0$ of Proposition $8$ satisfy to this
%property, except the two self-conjugate Specht modules
%$S^{(0),(3,2,1)}$ and $S^{(0),(4,1^3)}$. We deal with these modules
%first.
Suppose there exists in $V_6$ an invariant subspace that is
isomorphic to $S^{(0),(3,2,1)}$. Since $16>10$, we know that $\W\cap
V_5\neq\lb 0\rb$. By Lemma $3$ point $(ii)$ of \cite{CL},
$S^{(0),(2,2,1)}$ cannot occur in $V_5$ and by Results $3$ or $4$,
$S^{(0),(3,1,1)}$ cannot occur in $V_5$, hence $\W\cap V_5\simeq
S^{(0),(3,2)}$. This implies $l=r^3$ by Lemma $3$, point $(i)$ of
\cite{CL}. But when $l=r^3$, we have $k(6)=9$. Then, it is
impossible to have an invariant subspace of dimension $16$. So,
$S^{(0),(3,2,1)}$ does not occur in $V_6$ and consequently $\h{6}$
holds. Next, since $20>12$, it is an immediate consequence of
Results $3$ and $4$ and the branching rule that $S^{(0),(4,1^3)}$
cannot occur in $V_7$.
%We claim that $l\in\lb
%r^3,-r^3\rb$. Indeed, we have seen that when $l=\unsurr$, there
%exists in the Cohen-Wales space $V_n$ an irreducible
%$\chl$-dimensional invariant subspace isomorphic to
%$S^{(1,1),(n-2)}$ for all $n\geq 4$. When summing this dimension and
%the dimension of a Specht module from the list, one exceeds the
%dimension of the Cohen-Wales space. Thus, we cannot have
%$l=\unsurr$.
If $S^{(0),(4,2,1)}$ or its conjugate occur, then it is
impossible to have $l=r^3$, as otherwise, with Mathematica,
$k(7)=14<35$. %Then $l=-r^3$.
So, if $\W\subset V_7$ is such that $\W\simeq S^{(0),(4,2,1)}$, it
follows that $\W\cap V_6\simeq S^{(0),(4,1,1)}$. This is forbidden
by Result $3$. Also, if $\W\simeq S^{(0),(3,2,1^2)}$, then $\W\cap
V_6\simeq S^{(0),(3,1^3)}$ or $\W\cap V_6\simeq S^{(0),(3,2,1)}$, as
$S^{(0),(2^2,1^2)}$ cannot occur in $V_6$ by Lemma $3$ point $(ii)$
of \cite{CL}. But $S^{(0),(3,1^3)}$ cannot occur in $V_6$ by Result
$4$ and we have seen above that $S^{(0),(3,2,1)}$ cannot occur in
$V_6$ either. So, we are done with the second column and we have
shown that $\h{7}$ holds at this point. We now deal with the sixth
column. Suppose $\W\subset V_8$ with $\W\simeq S^{(0),(5,1^3)}$.
Then, by using Result $3$, we must have $\W\cap V_7\simeq
S^{(0),(4,1^3)}$. But we have seen above that $S^{(0),(4,1^3)}$
cannot occur in $V_7$. Similarly, by using Result $4$, the Specht
module $S^{(0),(4,1^4)}$ cannot occur in $V_8$. Hence $\h{8}$ holds.
Then, by the proof of Theorem $4$, if $l=r^3$, we have $k(8)=20$.
Suppose that $\W\subset V_9$ is isomorphic to one of the Specht
modules present in columns $3$ and $7$. Then, $k(9)\geq 48$ and so
$$\text{dim}(K(9)\cap V_8)\geq 48-16=32>20$$
Thus, it is impossible to have $l=r^3$.  %and so $l=-r^3$.
We use this fact to rule out the first Specht module present in
column $3$. For column $7$, we conclude directly as follows. If the
Specht modules of column $6$ cannot occur in $V_8$, then those of
column $7$ cannot occur in $V_9$ by using Results $3$ and $4$. By
the same arguments, the last column vanishes in turn. We now deal
with column $3$. Since $l\neq r^3$, the Specht module
$S^{(0),(6,2)}$ cannot occur in $V_8$ by Lemma $3$ point $(i)$ of
\cite{CL}. Since $S^{(0),(5,3)}$ can also not occur in $V_8$ by the
proof of Theorem $2$, we see that $S^{(0),(6,3)}$ cannot occur in
$V_9$. Next, if $\W\simeq S^{(0),(2^3,1^3)}$, by using the proof of
Theorem $4$ of \cite{CL} and Lemma $3$ point $(ii)$ of \cite{CL}, we
see that $\W\cap V_7\simeq S^{(0),(2^3,1)}$. Further, we have
$48>42=72-30$, so that $\W\cap V_6\neq \lb 0\rb$. By using again
Lemma $3$ point $(ii)$ of \cite{CL}, we get $\W\cap V_6\simeq
S^{(0),(2^3)}$. But this is impossible as shown in \cite{CL} at the
end of $\S\,3.4$. At this stage, having successfully ruled out
columns $3$ and $7$, we conclude that $\h{9}$ holds. It remains to
deal with column $4$. By the fact that $\h{9}$ holds, if $l=r^3$, we
have $k(9)=27$, as in the proof of Theorem $4$. If an invariant
subspace of dimension $75$ occurs in $V_{10}$, its intersection with
$V_9$ has a dimension greater than or equal to $57$ and must be
contained in $K(9)$. Therefore, it is impossible to have $l=r^3$.
This settles $S^{(0),(7,3)}$ by column $3$ and Lemma $3$, point
$(i)$ of \cite{CL}. As for $S^{(0),(2^3,1^4)}$, the conclusion is
even more straightforward and follows directly from column $3$ and
Lemma $3$,
point $(ii)$ of \cite{CL}.\\
We now process the other rows. Using Proposition $8$ and our work on
the row $k=0$, we keep in mind that Theorems $4$ and $6$ now hold,
always. Suppose $\W\simeq S^{(1),(3,3)}$. By Theorem $6$, the
existence of an irreducible invariant subspace of dimension $35$
forces $l=-r^3$. Moreover, since $S^{(0),(3,3)}$ cannot occur in
$V_6$ as seen several times in the past, we have $\W\cap V_6\simeq
S^{(1),(3,2)}$. Now the contradiction comes from
$$(\W\cap V_6)\da_{\mcalh(D_5)}\simeq S^{(0),(3,2)}\oplus S^{(1),(2,2)}\oplus S^{(1),(3,1)}$$
Indeed, by the proof of Lemma $3$, point $(i)$, the presence of
$S^{(0),(3,2)}$ in the restriction module $\W\cap
V_6\da_{\mcalh(D_5)}$ forces $l=r^3$. For $S^{(1),(2^3)}$, the proof
is similar. We use the fact that $S^{(0),(2^3)}$ does not occur in
$V_6$ (see the end of $\S\,3.4$ of \cite{CL}) and the fact that
$S^{(0),(2,2,1)}$ cannot be a constituent of $(\W\cap
V_6)\da_{\mcalh(D_5)}$ by the proof of Lemma $3$, point $(ii)$ of
\cite{CL}. So, we are done with the row $k=1$.\\
We now deal with the first column of the row $k=3$. Let $\W\subset
V_6$ such that $\W\simeq S^{(3),(1^3)}$. Then, by the branching
rule,
$$\W\da_{\mcalh(D_4)}\simeq S^{(1),(1^3)}\oplus
2\,S^{(2),(1^2)}\oplus S^{(1),(3)}$$ This is impossible by the end
of $\S\,3.5$ of \cite{CL}, where it is proven that $S^{(1),(1^3)}$
cannot occur in any Cohen-Wales space $V_n$. We now deal with
columns $2$ and $3$. Suppose first $\W\subset V_6$ is such that
%$\W\simeq S^{(3,3)^{+}}$ or $\W\simeq S^{(3,3)^{-}}$. We have
%$$S^{(3,3)^{+}}\da_{\mcalh(D_5)}\simeq
%S^{(3,3)^{-}}\da_{\mcalh(D_5)}\simeq S^{(3),(2)}$$ and
%$$S^{(3),(2)}\da_{\mcalh(D_4)}\simeq S^{(2,2)^{+}}\oplus
%S^{(2,2)^{-}}\oplus S^{(1),(3)}$$ Going back to the proof of Theorem
%$5$ of \cite{CL}, we see that this forces $l=-\unsur{r^7}$. By
%Theorem $3$ of the current paper, we then get $k(6)=6<10$. This
%contradicts $\W\subseteq K(6)$ and $\text{dim}(\W)=10$. If now
$\W\simeq S^{(1^3,1^3)^{+}}$ or $\W\simeq S^{(1^3,1^3)^{-}}$. We get
$$\W\da_{\mcalh(D_4)}\simeq S^{(1^2,1^2)^{+}}\oplus
S^{(1^2,1^2)^{-}}\oplus S^{(1),(1^3)}$$ By going back to the proof
of Theorem $5$ of \cite{CL}, we see that this is impossible:
$S^{(1),(1^3)}$ cannot be a constituent of $\W\da_{\mcalh(D_4)}$.
Note this same argument rules out $S^{(1^3),(1^4)}$ in column $4$.
To show that $S^{(3,3)^{+}}$ and $S^{(3,3)^{-}}$ don't occur in
$V_6$, we use computer means. Suppose at least one of them occurs in
$V_6$. Since $k(6)\geq 10$, then by using Mathematica we must have
$k(6)=15$ and $l=\unsurr$ or $k(6)\geq 24$ and $l=-r^3$. If the
first situation occurs, then there must exist an irreducible
$5$-dimensional invariant subspace in $V_6$, which also forces
$l=\unsur{r^5}$ by Theorem $4$ of \cite{CL}. Hence a contradiction
in that case. The dimensions of the irreducible
$\mcalh(D_6)$-modules that may occur in $V_6$ are
$$1,5,6,9,10,15,24$$
Suppose $l=-r^3$. For these values of $l$ and $r$, there cannot
exist any irreducible invariant subspace of dimension $5$ (resp $6$,
resp $9$, resp $15$), as otherwise $l=\unsur{r^5}$ (resp
$l=-\unsur{r^7}$, resp $l=r^3$, resp $l=\unsurr$) by Theorem $4$ of
\cite{CL} (resp Theorem $5$ of \cite{CL}, resp Theorem $6$ of
\cite{CL}, resp Theorem $2$ of the current paper). Nor can there
exist an irreducible invariant subspace of $V_6$ of dimension $24$
since $10+24=34>30=6\times 5$. Then it is impossible to have
$k(6)\geq 24$. We are now ready to process the rest of the fourth
column. Suppose $\W\subset V_7$ is such that $\W\simeq S^{(3),(4)}$.
Then we have
$$\W\da_{\mcalh(D_6)}\simeq S^{(3,3)^{+}}\oplus S^{(3,3)^{-}}\oplus
S^{(2),(4)}$$ On one hand, since $35>12$, we get that $\W\cap V_6$
is isomorphic to $S^{(2),(4)}$. Then, the existence of an
irreducible $15$-dimensional invariant subspace in $V_6$ forces
$l=\unsurr$ by Theorem $2$. On the other hand, the fact that $\W$
has dimension $35$ forces $l=-r^3$ by Theorem $6$. We obtain a
contradiction. So we are done with column $4$. Further, it is easy
to rule out $S^{(1^3),(4)}$ (resp $S^{(3),(1^4)}$) in column $5$ by
noticing that $35>12$ and by using the first column and Theorem $2$
(resp Result $2$) with $n=6$. Finally, by using the fourth column of
the row $k=3$ and the fact that $35>14$, we see that none of the
Specht modules present in the last row can occur in $V_8$. This
finishes the proof of Theorem $1$.

\end{document}